\documentclass[12pt,a4paper]{article}
\usepackage{latexsym,amsfonts,graphicx}

\newtheorem{thm}{Theorem}[section]
\newtheorem{lem}{Lemma}[section]
\newtheorem{cor}{Corollary}[section]
\newtheorem{prop}{Proposition}[section]

\newtheorem{hypo}{Hypothesis}[section]

\newenvironment{proof}{\vspace{-2ex} \trivlist \item[\hskip
\labelsep{\em Proof: }]}{{\phantom{\,}\hfill$\Box$}\endtrivlist}

\newcommand{\cpy}{\mbox{cap}}
\newcommand{\m}{\omega}
\newcommand{\E}{{\cal E}}
\newcommand{\F}{{\cal F}}
\newcommand{\K}{{\cal K}}
\newcommand{\LB}{\Delta}
\newcommand{\M}{{\cal M}}
\newcommand{\U}{{\cal U}}

\newcommand{\R}{\mathbb{R}}

\title{Discreteness of the spectrum of the Laplace-Beltrami operator}
\author{Mark Harmer\\
Department of Mathematics\\
Australian National University\\
Australia\\
email: mark.harmer@maths.anu.edu.au}
\date{}

\begin{document}

\maketitle 

\begin{abstract}
We propose simple conditions equivalent to the discreteness of the spectrum of the Laplace-Beltrami operator on a class of Riemannian manifolds close to warped products. For this class of manifolds we establish a relationship between discreteness of the spectrum and stochastic incompleteness.
\end{abstract}

\section{Introduction}
We say that an operator has only discrete spectrum if the essential spectrum of the operator is empty (i.e. the spectrum is discrete of finite multiplicity with no accumulation points). Much work has been done on finding conditions for the discreteness of the spectrum of the Laplace-Beltrami operator on functions, see for instance \cite{Don:Li, Bai, Kle, Bru, Kle2, Kon:Shu, KMS, Maz:Shu, Maz} and references therein.  Some of the most general necessary and sufficient conditions were found by Baider \cite{Bai} who considers manifolds which are warped products with the discreteness condition given as a condition on the spectrum of an operator defined on one of the terms in the product. Kleine and Br\"{u}ning \cite{Kle, Bru} express the Laplacian as a Schr\"{o}dinger operator with operator valued potential and derive general discreteness conditions for such operators. In \cite{Kle2} multiply warped products are considered and various necessary and sufficient conditions are derived in terms of the capacity. Very general conditions for discreteness are given by \cite{Kon:Shu, KMS, Maz:Shu, Maz} in terms of capacity. Our condition for discreteness is related to a condition on the capacity of subsets and may be considered as a particular case of Maz'ja's condition \cite{Maz, Grig2}. \\
For us an important condition for discreteness was given by Kac and Krein \cite{KaKr} for the singular string. The main result in our paper may be considered as a generalisation of the result of Kac and Krein to manifolds with ends which are `close to' warped products. \\
In this paper we prove a condition equivalent to discreteness of the spectrum of the Laplace-Beltrami operator on a manifolds with ends having a certain type of bound on the mean curvature. Analogous to the result of Kac and Krein our condition splits into two cases depending on whether the end of the manifold has finite or infinite volume. We use ideas which appear in \cite{KaKr} and a paper by Muckenhoupt \cite{Muc} on weighted Hardy inequalities. The conditions equivalent to discreteness  allow one to easily see the relationship between volume growth down the ends of the manifold and discreteness. \\
We follow the main result with a discussion on the relationship between our condition and the condition of Maz'ja based on capacity \cite{Maz}. The paper concludes by showing that, for the chosen class of manifolds, stochastic incompleteness is sufficient for discreteness. 
\section{Preliminaries}
We consider a complete, non compact $n+1$-dimensional Riemannian manifold $\M$. We impose two hypotheses on $\M$.
\begin{hypo}\label{hyp1}
There is an compact subset $\U$ of $\M$ such that $\M\setminus\U$ consists of a finite number of disjoint, noncompact ends $\E_i$. Each end $\E_i$ is diffeomorphic to $\R_+\times\K_i$ where $\K_i$ is an $n$-dimensional compact manifold. The diffeomorphism induces on $\R_+ \times \K_i$ the metric
\begin{equation}\label{bass1}
ds^2 = dr^2 + d\theta_{\K_i} (r)^2  
\end{equation}
where $d\theta_{\K_i} (r)^2$ is a smooth family of metrics on $\K_i$. 
\end{hypo}
The decomposition principle (see Proposition 2.1 of \cite{Don:Li}, Lemma 2.3 of \cite{Bai} or \cite{Gla} for a precise statement) states that the essential spectrum of the Laplacian on functions is invariant under compact perturbation. Consequently, the discreteness of the Laplacian on $\M$ is equivalent to the discreteness of the Dirichlet Laplacian on each of the ends. For this reason we drop the subscript $i$ and consider questions of discreteness on a generic end $\E$ in the sequel. \\
The form of volume on $\E$ may be written 
$$
d\mu = \m (r,\theta)\, d\theta^{n}\, dr
$$
while on $\E$ the Laplace-Beltrami operator on functions has the form
$$
\LB \equiv - \frac{1}{\m} \frac{\partial}{\partial r} \m \frac{\partial}{\partial r} + \Delta_{\K} (r)
$$
where $\Delta_{\K} (r)$ is the Laplacian on $\K$ with respect to the metric $d\theta_{\K} (r)^2$. We denote the norm in $L_2 \left( \M , d\mu \right)$ by $\| \cdot \|$.\\
The quantity
$$
h (r, \theta) = \frac{\m'}{\m} (r, \theta)
$$
is the mean curvature of $\K$ at the point $(r, \theta)$ \cite{Gage}. In this paper $\mbox{}'$ will denote differentiation with respect to $r$. This leads us to our second hypothesis:
\begin{hypo}\label{hyp2}
The mean curvature $h$ of $\K$ satisfies
\begin{equation}\label{bass2}
\sup_{\E} \left| h - \frac{1}{\bar{\m}} \int_{\K} h \m \, d\theta^n \right| = \sup_{\E} \left| h - \bar{h} \right| < c
\end{equation}
where $\bar{\m} (r) = \int_{\K} \m \, d\theta^n$ is the volume of $\K$ at the point $r$ and $\bar{h}$ is defined in the equality. 
\end{hypo}
This assumption means that the end is very close to a warped product. This allows us to effectively reduce the problem to a one-dimensional one: the end $\E$ is sufficiently well behaved that only the $r$ dependence of the metric is important for discreteness. The mean curvature is not an intrinsically defined quantity so it is not clear whether hypothesis \ref{hyp2} is intrinsic or not. \\
We define the non compact subset $\E_t \subset \E$, $t\ge0$, as the set diffeomorphic to
$$
(t,\infty) \times \K  
$$
under the diffeomorphism defined above and we denote the norm in $L_2 \left( \E_t , d\mu \right)$ by $\| \cdot \|_t$. \\
Given $f\in C^{\infty}_0 (\E,\R)$ we define (up to sign) the following average
\begin{equation}\label{ave}
\bar{f} (r) = \left( \frac{1}{\bar{\m} (r)} \int_{\K} f^2  \m  \, d\theta^n \right)^{\frac{1}{2}} \, .
\end{equation}
We say that such an averaged function is spherically symmetric. 
\begin{lem}\label{blem}
Given $f\in C^{\infty}_0 (\E_t ,\R)$ on an end $\E$ satisfying hypothesis \ref{hyp2} we have
\begin{equation}\label{coer}
\frac{1}{2} \left\| \bar{f}' \right\|^2_{t} - \frac{c^2}{4} \left\| f \right\|^2_{t} \le \left\| \nabla f \right\|^2_{t} \, .
\end{equation}
\end{lem}
\begin{proof}
We differentiate (\ref{ave}) 
\begin{eqnarray*}
\left| \bar{f}' \right| & \le & \frac{1}{2\bar{f}} \left[  \left| \frac{2}{\bar{\m}} \int_{\K} f f' \m \, d\theta^n \right| + \left| \frac{1}{\bar{\m}} \int_{\K} f^2 \m' \, d\theta^n - \bar{h} \bar{f}^2 \right| \right] \\
& \le & \frac{1}{2\bar{f}} \left[  \frac{2}{\bar{\m}} \left[ \int_{\K} f^2 \m \, d\theta^n \right]^{\frac{1}{2}} \left[ \int_{\K} \left(f'\right)^2 \m \, d\theta^n \right]^{\frac{1}{2}} + \sup_{\K} \left| h - \bar{h} \right| \bar{f}^2 \right] \\
& \le & \left[ \frac{1}{\bar{\m}} \int_{\K} \left(f'\right)^2 \m \, d\theta^n \right]^{\frac{1}{2}}
+ \frac{c}{2} \bar{f} \, .
\end{eqnarray*}
Using $(a+b)^2\le 2a^2+2b^2$ we get 
$$
\frac{1}{2} \left(\bar{f}'\right)^2 \bar{\m} - \frac{c^2}{4} \bar{f}^2 \bar{\m} \le \int_{\K} \left(f'\right)^2 \m \, d\theta^n \, .
$$
Then integrating over $r$ and using 
$$
 \int \int_{\K} \left(f'\right)^2 \m \, d\theta^n \, dr \le \int \left[ \int_{\K} \left(f'\right)^2 \m \, d\theta^n + \int_{\K} \left| \nabla_{\K} f \right|^2 \m \, d\theta^n \right] \, dr = \int_{\E} \left| \nabla f \right|^2 \, d\mu \, ,
$$
where $\nabla_{\K}$ is the gradient on $\K$, we get the result.
\end{proof}
For discreteness the stronger condition (\ref{bass2}) is not necessary, all we need is the condition (\ref{coer}) of the lemma, indeed we may choose this as our second hypothesis. We refer to (\ref{coer}) as the `coerciveness' condition (compare with equation (1.13) of \cite{Bru}). Ideally though we would replace hypothesis \ref{hyp2} with an intrinsic condition. \\
We define the Dirichlet \cite{Ura}
$$
\lambda_0 \left( \E_t \right) = \inf \left\{ \frac{\left\| \nabla f \right\|^2_{t}} 
{\left\| f \right\|^2_{t}} \, : \: f,\LB f \in L_2 \left( \E_t \right) \, , \: f \in H^1_0 \left( \E_t \right) \right\} 
$$
and Neumann
$$
\mu_0 \left( \E_t \right) = \inf \left\{  \frac{\left\| \nabla f \right\|^2_{t}} 
{\left\| f \right\|^2_{t}} \, : \: f,\LB f \in L_2 \left( \E_t \right) \, , \: \left( \LB f , u \right) = \left( \nabla f , \nabla u \right) \: \forall u \in H^1 \left( \E_t \right) \right\} 
$$
Rayleigh quotients where $H^1 \left( \E_t \right)$ is the set of functions with $\| u \|^2_t + \| \nabla u \|^2_t < \infty$ and $H^1_0 \left( \E_t \right)$ is the closure of $C^{\infty}_0 \left( \E_t \right)$ in $H^1 \left( \E_t \right)$.

The following result (\cite{Bai}, theorem 2.2) is our main tool to prove discreteness:
\begin{thm}\label{Bai1}
The spectrum of the Laplace-Beltrami operator $\LB$ will be discrete iff on each end $\E$ the Dirichlet Rayleigh quotient satisfies
$$
\lim_{t\to\infty} \lambda_0 \left( \E_t \right) = +\infty \, . 
$$
\end{thm}
It is easy to see that the same result will not hold for the Neumann Rayleigh quotient (consider a finite volume non compact end down which the Dirichlet Rayleigh quotient is unbounded but where the Neumann Rayleight quotient remains bounded because of the presence of the constant eigenfunction). Nevertheless, in the infinite volume case we have:
\begin{lem}\label{Bai2}
The spectrum of the Laplacian $\LB$ on a manifold $\M$ satisfying hypotheses \ref{hyp1} and \ref{hyp2} will be discrete iff on each finite volume end 
\begin{equation}\label{dcnd}
\lim_{t\to\infty} \lambda_0 \left( \E_t \right) = +\infty 
\end{equation}
and on each infinite volume end
\begin{equation}\label{ncnd}
\lim_{t\to\infty} \mu_0 \left( \E_t \right) = +\infty \, .
\end{equation}
\end{lem}
\begin{proof}
We need to establish that on an infinite volume end the unboundedness of the Neumann quotient is equivalent to the unboundedness of the Dirichlet quotient. The Dirichlet quotient is bounded below by the Neumann quotient so that unboundedness of the Neumann quotient implies unboundedness of the Dirichlet quotient. \\
Conversely, suppose that there exists a $\mu<\infty$ such that for all $t>0$, $\mu_0 \left( \E_t \right) < \mu$. This means we can find a $\psi_1$ satisfying the above conditions such that
$$
\frac{ \| \nabla \psi_1 \|^2_0 }{ \| \psi_1 \|^2_0 } < \mu \, .
$$
Furthermore, since $\E_t$ has infinite volume, we can choose $\psi_1$ so that it has compact support in $r\in[0,t_{1})$ for some $t_1$. Repeating this at $t = t_1$  we get a sequence $\{ \psi_i , t_i \}$ such that $\psi_i$ satisfies the above inequality and has support in $r\in[t_{i-1},t_{i})$. We average the $\psi_i$ as in (\ref{ave}) to get a sequence of spherically symmetrical functions $\bar{\psi}_i$ where, without loss of generality we may assume that
$$
\left. \bar{\psi}_i \right|_{r=t_{i-1}} = 1 \, , \quad \left. \bar{\psi}'_i \right|_{r=t_{i-1}} = 0 \, ,
$$
and we continue each $\bar{\psi}_i$ in $r<t_{i-1}$ as $\bar{\psi}_i =1$. In forming the average (defined up to signs) of a complex valued function, we take the average of the real and imaginary parts separately. \\
Clearly we can take a subsequence $\{ \bar{\psi}_{k_l} , t_{k_l} \}$ such that
$$
\left\| \bar{\psi}_{k_{l+1}} \right\|_{t_{k_l}} > \left\| \bar{\psi}_{k_{l}} \right\|_{t_{k_{l-1}}} \, .
$$
We drop the extra subscript and denote this subsequence by the same notation $\{ \psi_l , t_l \}$ so that 
$$
\left\| \bar{\psi}_{{l+1}} \right\|_{t_{l}} > \left\| \bar{\psi}_{{l}} \right\|_{t_{{l-1}}} \, .
$$
Furthermore, it is clear from the definition and (\ref{coer}) that 
$$
\frac{\left\| \nabla \bar{\psi}_l \right\|^2_{t_{l-1}} }{ \left\| \bar{\psi}_{l} \right\|^2_{t_{l-1}} } < 2\mu + \frac{c^2}{2} \, .
$$
We consider the sequence $\varphi_l = \bar{\psi}_{2l} - \bar{\psi}_{2l-1}$ defined on the whole end $\E$. This sequence has disjoint supports and
\begin{eqnarray*}
\frac{\left\| \nabla \varphi_l \right\|^2_{0} }{ \left\| \varphi_{l} \right\|^2_{0} } & = & \frac{\left\| \nabla \bar{\psi}_{2l} \right\|^2_{t_{2l-1}} + \left\| \nabla \bar{\psi}_{2l-1} \right\|^2_{t_{2l-2}}}{ \left\| \varphi_{l} \right\|^2_{0} } \\
& \le & \frac{ \left\| \nabla \bar{\psi}_{2l} \right\|^2_{t_{2l-1}} }{ \left\| \bar{\psi}_{2l}  \right\|^2_{t_{2l-1}} } + \frac{ \left\| \nabla \bar{\psi}_{2l-1} \right\|^2_{t_{2l-2}} }{ \left\| \bar{\psi}_{2l-1}  \right\|^2_{t_{2l-2}} } \le 4\mu + c^2 \, .
\end{eqnarray*}
In the first line we use the fact that $\nabla \bar{\psi}_{2l}$ and $\nabla \bar{\psi}_{2l-1}$ have disjoint supports while in the second line we use $\left\| \varphi_{l} \right\|_{0} > \left\| \bar{\psi}_{{2l}} \right\|_{t_{2l-1}} > \left\| \bar{\psi}_{{2l-1}} \right\|_{t_{{2l-2}}}$. Consequently, according to \cite{Gla}, theorem 13 page 15, we have a point of essential spectrum.
\end{proof}
We require one more result which  will be used in the discussion of brownian motion.
\begin{lem}\label{pertlem}
Suppose that $\M$ is a manifold satisfying hypothesis \ref{hyp1} with only one end and such that the Laplacian has non empty essential spectrum. In particular the metric down the end is of the form (\ref{bass1}). Perturbing the metric by an exponential factor to
\begin{equation}\label{pertm}
ds^2 = dr^2 + e^{2cr/n} d\theta_{\K} (r)^2 \, ,
\end{equation}
where $c>0$, the essential spectrum of the Laplacian will remain non empty.
\end{lem}
\begin{proof}
The assumption of non empty essential spectrum implies, by theorem \ref{Bai1}, the existence of $\lambda<\infty$ and a family $\varphi_t\in L_2 \left( \E_t \right)$ such that
$$
\frac{\left\| \nabla \varphi_t \right\|^2_{t} }{ \left\| \varphi_{t} \right\|^2_{t} } < \lambda \, .
$$
We define
$$
\varphi_{t,c} = e^{-cr/2} \varphi_{t}
$$
and denote the norm associated to the perturbed metric (\ref{pertm}) by $\| \cdot \|_{t,c}$. Then $\| \varphi_{t,c} \|_{t,c} = \| \varphi_{t} \|_{t}$ and
\begin{eqnarray*}
\left\| \nabla_c \varphi_{t,c} \right\|^2_{t,c} & = & \int_{\E_t} \left( \varphi'_{t,c} \right)^2 e^{cr} \m \, d\theta^n \, dr +  \int_{\E_t} \left| \nabla_{\K} \varphi_{t,c} \right|^2 e^{-2cr/n} e^{cr} \m \, d\theta^n \, dr \\
& = & \int_{\E_t} \left( \varphi'_{t} - \frac{c}{2} \varphi_{t} \right)^2 \m \, d\theta^n \, dr +  \int_{\E_t} \left| \nabla_{\K} \varphi_{t} \right|^2 e^{-2cr/n} \m \, d\theta^n \, dr \\
& \le & 2 \left\| \nabla \varphi_t \right\|^2_t + \frac{c^2}{2} \left\| \varphi_t \right\|^2_t
\end{eqnarray*}
which implies that
$$
\frac{\left\| \nabla_c \varphi_{t,c} \right\|^2_{t,c}}{\| \varphi_{t,c} \|^2_{t,c}} \le 2 \lambda + \frac{c^2}{2} \, .
$$
\end{proof}
\section{Main result}
\begin{thm}\label{mthm}
The spectrum of the Laplacian on a manifold $\M$ satisfying hypotheses \ref{hyp1} and \ref{hyp2} is discrete iff for each end $\E$ either
\begin{equation}\label{ldc1}
\lim_{t\to\infty} \sup_{s>t} \int^{s}_t \bar{\m}^{-1} \, dr \int^{\infty}_s \bar{\m} \, dr = 0 \, , 
\end{equation} 
or 
\begin{equation}\label{ldc2}
\lim_{t\to\infty} \sup_{s>t} \int^{\infty}_s \bar{\m}^{-1} \, dr \int^{s}_t \bar{\m} \, dr = 0 \, . 
\end{equation}
\end{thm}
\begin{proof} 
The proof is split into two cases: either $\int^{\infty}_0 \bar{\m}^{-1} \, dr = \infty$ and we show that (\ref{ldc1}) is equivalent to the unboundedness of the Dirichlet Rayleigh quotient; or $\int^{\infty}_0 \bar{\m} \, dr = \infty$ and we show that (\ref{ldc2}) is equivalent to the unboundedness of the Neumann Rayleight quotient. We refer to these as the finite and infinite volume cases respectively.\\
\begin{enumerate} 
\item Assuming $\int^{\infty}_0 \bar{\m}^{-1} \, dr = \infty$ we show that (\ref{ldc1}) is equivalent to (\ref{dcnd}). \\
We choose $u \in C^{\infty}_0 \left( \E_t \right)$ assuming, without loss of generality, that $u$ is real. 
Then
\begin{eqnarray*}
\int_{\E_t} u^2 \, d\mu & = & \int^{\infty}_{t} \bar{u}^2 \, \bar{\m}\, dr = 2 \int^{\infty}_{t} \int^{\infty}_{r} \bar{\m} \, ds \, \bar{u} \bar{u}' \, dr \\
& = & 2 \sup_{s>t} \left( \int^{s}_t \bar{\m}^{-1} \, dr \int^{\infty}_s \bar{\m} \, dr \right)
\int^{\infty}_{t} \frac{1}{\int^{r}_{t} \bar{\m}^{-1} \, ds} \bar{u} \bar{u}' \, \, dr \\
& \le & 2 \sup_{s>t} \left( \int^{s}_t \bar{\m}^{-1} \, dr \int^{\infty}_s \bar{\m} \, dr \right) \times \mbox{} \\
& & \left[ \int^{\infty}_{t} \bar{u}^2 \, \frac{\nu^2 (r,t)}{\bar{\m}} \, dr \right]^{\frac{1}{2}}
\left[ \int^{\infty}_{t} \left(\bar{u}' \right)^2 \, \bar{\m} \, dr \right]^{\frac{1}{2}} \, . 
\end{eqnarray*}
Here we have denoted
\begin{eqnarray*}
\nu(r,t) & = & \left( \int^{r}_{t} \bar{\m}^{-1} \, ds \right)^{-1} \\
\Rightarrow \nu' (r,t) & = & -\frac{\nu^2 (r,t)}{\bar{\m} (r)} \\
\Rightarrow \nu (r,t) & = & \int^{\infty}_r \frac{\nu^2 (s,t)}{\bar{\m} (s)} \, ds \, .
\end{eqnarray*}
Using this in the first integral on the right hand side of our inequality we have
\begin{eqnarray*}
\int^{\infty}_{t} \bar{u}^2 \, \frac{\nu^2 (r,t)}{\bar{\m}} \, dr & = & 2 \int^{\infty}_{t} \bar{u} \bar{u}' \, \nu (r,t) \, dr \\
& \le & 2 \left[ \int^{\infty}_{t} \bar{u}^2 \, \frac{\nu^2 (r,t)}{\bar{\m}} \, dr \right]^{\frac{1}{2}}
\left[ \int^{\infty}_{t} \left(\bar{u}'\right)^2 \, \bar{\m} \, dr \right]^{\frac{1}{2}} \, , 
\end{eqnarray*}
that is
$$
\int^{\infty}_{t} \bar{u}^2 \, \frac{\nu^2 (r,t)}{\bar{\m}} \, dr \le 4 \int^{\infty}_{t} \left(\bar{u}'\right)^2 \, \bar{\m} \, dr \, .
$$
Putting this back into the original inequality we have
$$
\int_{\E_t} u^2 \, d\mu \le 4 \sup_{s>t} \left(\int^{s}_t \bar{\m}^{-1} \, dr \int^{\infty}_s \bar{\m} \, dr \right) 
\int^{\infty}_{t} \left(\bar{u}'\right)^2 \, \bar{\m} \, dr \, . 
$$
On the other hand, using lemma \ref{blem} we can write the Dirichlet integral as
$$
\frac{1}{2} \int \left(\bar{u}'\right)^2 \bar{\m}\, dr - \frac{c^2}{4} \int \bar{u}^2 \bar{\m} \, dr \le \int \left| \nabla u \right|^2 \, d\mu \, .
$$
Putting these inequalities together we have
$$
\left( \frac{1}{8} \frac{1}{\sup_{s>t} \left( \int^{s}_t \bar{\m}^{-1} \, dr \int^{\infty}_s \bar{\m} \, dr \right)} - \frac{c^2}{4} \right) \int_{\E_t } {u}^2 \, d\mu \le \int_{\E_t } \left| \nabla u \right|^2 \, d\mu
$$
which proves that (\ref{ldc1}) is sufficient for discreteness. \\
For necessity we consider the family $ v(r;t,s,s_0,s_1) $, $t<s<s_0<s_1$, of Lipschitz functions
$$
v (r) = \left\{ \begin{array}{cl}
\frac{\int^r_t \bar{\m}^{-1} \, dr'}{\int^s_t \bar{\m}^{-1} \, dr'} & : t<r<s \\
1 & : s<r<s_0 \\
2 - \left(\frac{\int^r_t \bar{\m}^{-1} \, dr' - \int^s_t \bar{\m}^{-1} \, dr'}{\int^{s_1}_{s_0} \bar{\m}^{-1} \, dr'} \right) & :s_0<r<s_1 \\
0 & :s_1<r 
\end{array} \right.
$$
where, since $\int^\infty_0 \bar{\m}^{-1} \, dr$, we are able to choose $s_1$ so that
$$
\int^{s_1}_{s_0} \bar{\m}^{-1} \, dr' = \int^{s_0}_{s} \bar{\m}^{-1} \, dr' \, .
$$
Putting this into the Rayleigh quotient we see that
$$
\frac{\left\| \nabla v \right\|^2_t} 
{\left\| v \right\|^2_t} \le \left( 1 + \frac{\int^s_t \bar{\m}^{-1} \, dr'}{\int^{s_0}_s \bar{\m}^{-1} \, dr'} \right) 
\left( \int^s_t \bar{\m}^{-1} \, dr' \int^{s_0}_s \bar{\m} \, dr'\right)^{-1} \, .
$$
Letting $s_0\rightarrow\infty$ gives the result.
\item Using the assumption of infinite volume $\int^{\infty}_{0} \bar{\m} \, dr = \infty$ and the same argument as above but with $\bar{\m}\leftrightarrow\bar{\m}^{-1}$ we see that for $v \in C^{\infty}_0 \left( \E_t \right)$
\begin{equation}\label{infvol}
\left[ \int^{\infty}_{t} \bar{v}^2 \, \bar{\m}^{-1} \, dr \right]^{\frac{1}{2}} \le 2 
\sup_{s>t} \left( \int^{\infty}_s \bar{\m}^{-1} \, dr \int^{s}_t \bar{\m} \, dr \right)
\left[ \int^{\infty}_{t} \left(\bar{v}'\right)^2 \, \bar{\m}^{-1} \, dr \right]^{\frac{1}{2}} \, .
\end{equation}
Now consider again a real function $u \in C^{\infty}_0 \left( \E_t \right)$;
\begin{eqnarray*}
\int_{\E_t} u^2 \, d\mu & = & \int^{\infty}_{t} \bar{u}^2 \, \bar{\m} \, dr =  \int^{\infty}_{t} \int^{r}_{t} \bar{u} (s) \bar{\m} (s)\, ds \, \bar{u}' (r)\, dr \\
& = &  \int^{\infty}_{t} \frac{1}{\bar{\m}^{1/2} } \int^{r}_{t} \bar{u} \bar{\m} \, ds \, \bar{u}' \, \bar{\m}^{1/2} \, dr \\
& \le &  \left[ \int^{\infty}_{t} \left(\int^{r}_{t} \bar{u} \bar{\m} \, ds\right)^2 \, \bar{\m}^{-1} \, dr \right]^{\frac{1}{2}} \left[ \int^{\infty}_{t} \left(\bar{u}'\right)^2 \, \bar{\m} \, dr \right]^{\frac{1}{2}} \, . 
\end{eqnarray*}
Using (\ref{infvol}) on the right hand side with $\bar{v}(r)=\int^r_t \bar{u} \bar{\m} ds$ we get the inequality
$$
\int_{\E_t} u^2 \, d\mu \le 4 \sup_{s>t} \left( \int^{\infty}_s \bar{\m}^{-1} \, dr \int^{s}_t \bar{\m} \, dr \right) 
\int^{\infty}_{t} \left(\bar{u}'\right)^2 \, \bar{\m} \, dr  
$$
analogous to the result in the finite volume case. Again using lemma \ref{blem} we get that (\ref{ldc2}) is sufficient for discreteness. \\
For necessity we consider the family $ v(r;t,s,s_0) $, $t<s<s_0$, of Lipschitz functions
$$
v (r) = \left\{ \begin{array}{cl}
1 & : t<r<s \\
1 - \left( \frac{\int^r_t \bar{\m}^{-1} \, dr' - \int^s_t \bar{\m}^{-1} \, dr'}{\int^{s_0}_s \bar{\m}^{-1} \, dr'} \right) & : s<r<s_0 \\
0 & : s_0<r
\end{array} \right.
$$
in the Rayleigh quotient:
$$
\frac{\left\| \nabla v \right\|^2_t}{\left\|  v \right\|^2_t} \le \frac{1}{\int^{s_0}_s \bar{\m}^{-1} \, dr' \int^{s}_t \bar{\m} \, dr'} \, .
$$
Again letting $s_0\rightarrow\infty$ gives the result.
\end{enumerate}
\end{proof}
The proof of this theorem uses ideas from \cite{Bir, Maz}, the proof of discreteness for the singular string by Kac and Krein \cite{KaKr} and Muckenhoupt's weighted Hardy inequality \cite{Muc}. The Kac and Krein result can be directly used to give a condition for discreteness in the warped product case but, in the form it is published, is not convenient to use in the more general context considered here. Muckenhoupt gives a result very similar to this theorem (in somewhat greater generality); however, he requires finiteness of the arguments of the limits in (\ref{ldc1}, \ref{ldc2}) which we are able to relax thereby simplifying the proof. \\
We note that the conditions (\ref{ldc1},\ref{ldc2}) for discreteness have some resemblence to the condition for discreteness given by Br\"{u}ning in the case of the warped product (theorem 3.4 of \cite{Bru}). \\
These conditions for discreteness can be further simplified.
\begin{cor}\label{mcor}
The spectrum of the Laplacian on a manifold $\M$ satisfying hypotheses \ref{hyp1} and \ref{hyp2} is discrete iff for each end $\E$ either
\begin{equation}\label{sdc1}
\lim_{s\to\infty} \int^{s}_0 \bar{\m}^{-1} \, dr \int^{\infty}_s \bar{\m} \, dr = 0 \, , 
\end{equation} 
or 
\begin{equation}\label{sdc2}
\lim_{s\to\infty} \int^{\infty}_s \bar{\m}^{-1} \, dr \int^{s}_0 \bar{\m} \, dr = 0 \, . 
\end{equation}
\end{cor}
\begin{proof}
We consider the equivalence of (\ref{ldc2}) to (\ref{sdc2}), the other case follows from a similar argument. \\
Assuming (\ref{ldc2}) we can, given $\epsilon>0$, find $T$ such that for $t\ge T$
$$
\sup_{s>t} \int^{\infty}_s \bar{\m}^{-1} \, dr \int^s_t \bar{\m} \, dr < \epsilon \, .
$$
On the other hand (\ref{ldc2}) implies $\int^{\infty}_0 \bar{\m}^{-1} \, dr<\infty$ so that we can find $S$ such that for $s\ge S$ 
$$
\int^{\infty}_s \bar{\m}^{-1} \, dr < \frac{\epsilon}{M}
$$
where $M = \int^T_0 \bar{\m} \, dr$. Consequently, for $s\ge\max (T,S)$ we have
$$
\int^{\infty}_s \bar{\m}^{-1} \, dr \int^{s}_0 \bar{\m} \, dr = \int^{\infty}_s \bar{\m}^{-1} \, dr \left( \int^{T}_0 \bar{\m} \, dr + \int^{s}_T \bar{\m} \, dr \right) < 2 \epsilon \, .
$$
Now we suppose that (\ref{sdc2}) holds from which we have, for all $s>0$,
$$
\lim_{t\to\infty} \int^{\infty}_{s+t} \bar{\m}^{-1} \, dr \int^{t}_0 \bar{\m} \, dr = 0 
$$
following from the monotonicity of $\int^{\infty}_{t} \bar{\m}^{-1} \, dr$. Consequently, for all $s>0$,
\begin{eqnarray*}
0 & = & \lim_{t\to\infty} \left( \int^{\infty}_{s+t} \bar{\m}^{-1} \, dr \int^{s+t}_0 \bar{\m} \, dr - \int^{\infty}_{s+t} \bar{\m}^{-1} \, dr \int^{t}_0 \bar{\m} \, dr\right) \\
& = & \lim_{t\to\infty} \int^{\infty}_{s+t} \bar{\m}^{-1} \, dr \int^{s+t}_t \bar{\m} \, dr \\
& = & \lim_{t\to\infty} \int^{\infty}_{s} \bar{\m}^{-1} \, dr \int^{s}_t \bar{\m} \, dr \, , \qquad \forall s>t \, ,
\end{eqnarray*}
in particular (\ref{ldc2}) holds.
\end{proof}
If we assume that the metric down the ends of the manifold is sufficiently `well behaved' then it is possible to greatly simplify conditions (\ref{ldc1},\ref{ldc2}) or (\ref{sdc1},\ref{sdc2}). Specifically, if we assume that the following limits exist (we will say that a limit going to $\pm\infty$ exists) then it is clear that (\ref{sdc1},\ref{sdc2}) are equivalent to 
$$
\lim_{s\to\infty} \frac{d}{ds} \ln \int^{\infty}_s \bar{\m} \, dr =  -\infty 
$$
or
$$
\lim_{s\to\infty} \frac{d}{ds} \ln \int^{\infty}_s \bar{\m}^{-1} \, dr = -\infty
$$
respectively. This is a simple consequence of L'Hopital's rule \cite{Tay}, in the case of (\ref{sdc2})
\begin{eqnarray*}
0 & = & \lim_{s\to\infty} \int^{\infty}_s \bar{\m}^{-1} \, dr \int^{s}_0 \bar{\m} \, dr \\
& = & \lim_{s\to\infty} \frac{\int^{s}_0 \bar{\m} \, dr}{\left( \int^{\infty}_s \bar{\m}^{-1} \, dr \right)^{-1}} \\
& = & \lim_{s\to\infty} \frac{\bar{\m} (s)}{\bar{\m}^{-1} (s) \left( \int^{\infty}_s \bar{\m}^{-1} \, dr \right)^{-2}} 
\end{eqnarray*}
where the third equality follows because we have assumed that the limit exists. \\
Hypothesis \ref{hyp2} is clearly important for sufficiency of the above condition for discreteness. In the case of necessity hypothesis \ref{hyp2} only appears in lemma \ref{Bai2}. Nevertheless, if we drop hypothesis \ref{hyp2} we have the following proposition.
\begin{prop}
The Laplacian on a manifold $\M$ satisfying hypothesis \ref{hyp1} has essential spectrum if on one of the ends
$$
\int^{\infty}_0 \bar{\m} \, dr=\infty=\int^{\infty}_0 \bar{\m}^{-1} \, dr \, . 
$$
\end{prop}
\begin{proof}
We use the hypothesis to construct a characteristic sequence. Introducing a new coordinate (the map between a warped product and a singular string)
$$
z(r) = \int^{r}_{0} \bar{\m}^{-1} \,ds \, ,
$$
$z\in [0,\infty )$, we define $\hat{r}$ to be the unique solution of $z(\hat{r}) = 2z(r)$. We claim that our assumption implies the existence of a sequence $\{ r_l \}$ such that 
\begin{equation}\label{ndisc}
\int^{r_l}_{0} \bar{\m}^{-1} \, {ds} \int^{\hat{r}_l}_{r_l} \bar{\m}\, ds \ge c_0 > 0 
\end{equation}
By contradiction we suppose that given $\epsilon >0$ we can find $R$ such that for all $r>R$
$$
\int^{r}_{0} \bar{\m}^{-1}\, dt \int^{\hat{r}}_{r} \bar{\m}\, dt < \epsilon \, . 
$$
Then
\begin{eqnarray*}
\int^{r}_{0} \bar{\m}^{-1}\, ds \int^{\infty}_{r} \bar{\m}\, dx & = & \int^{r}_{0} \bar{\m}^{-1}\, ds \left( \int^{\hat{r}}_{r} \bar{\m}\, ds + \int^{\hat{\hat{r}}}_{\hat{r}} \bar{\m}\, ds + \cdots \right) \\
& = & z (r) \int^{\hat{r}}_{r} \bar{\m}\, ds + \frac{1}{2} 2 z(r) \int^{\hat{\hat{r}}}_{\hat{r}} \bar{\m}\, ds + \cdots \\
& \le & \epsilon + \frac{1}{2} \epsilon + \cdots + \frac{1}{2^n} \epsilon + \cdots = 2 \epsilon 
\end{eqnarray*}
which contradicts the assumption of infinite measure. \\
Let us take a subsequence $\{ r_{l_{k}} \}$ such that
\begin{equation}\label{orth}
z \left( r_{l_{k+1}} \right) > 3\, z \left( r_{l_{k}} \right) \, .
\end{equation}
We drop the extra subscript, denoting this sequence by $\{ r_k \}$, and also put $z_k = z(r_k)$.
Let us choose a smooth real function $\eta$ with support $(\frac{4}{5},\frac{11}{5})$ which is equal to one on $[1,2]$. Then we construct the sequence of functions 
$$
f_k (r, \theta) = \sqrt{z_k}\; \eta \left( \frac{z(r)}{z_k}\right) \, .
$$
The Dirichlet integral is 
$$
\int | f'_k (r)|^2 \bar{\m} \, dr = \int \left| \eta' \left(\frac{z (r)}{z_k} \right) \right|^2 \frac{1}{z_k} \bar{\m}^{-1}\, dr = \int |\eta' (y)|^2 dy = c_1
$$
for some constant $c_1$. Here $\eta'$ denotes differentiation with respect to the argument of $\eta$, not $r$. Furthermore, the $L_2$ norms of the $u_k$
\begin{eqnarray*}
\int |f_k|^2 \bar{\m} \, dr & = & z_k \int \left| \eta \left( \frac{z(r)}{z_k}\right) \right|^2 \bar{\m} \, dr
\ge z_k \int^{\hat{r}_k}_{r_k} \bar{\m} \, dr \\
& = & \int^{r_k}_{0} \bar{\m}^{-1}\, ds \int^{\hat{r}_k}_{r_k} \bar{\m}\, ds \ge c_0
\end{eqnarray*}
are bounded below. Therefore we may normalise this sequence in $H_x$, $g_k = f_k
/ \| f_k \|$, and in doing so the Dirichlet integral remains bounded by
$$
\int_{\E } \left| \nabla g_k \right|^2 \, d\mu \le \frac{c_1}{c_0} \, .
$$
Furthermore we have from (\ref{orth}) that the support of $g_l$, $g_m$, $l\ne m$, are disjoint so the sequence is orthonormal. Using the same argument as above we have a characteristic sequence and therefore a point of essential spectrum.
\end{proof}
Our conditions (\ref{ldc1},\ref{ldc2}) equivalent to discreteness can be thought of as a simplification of a general result of Maz'ja. For a compact subset $\F\subset\subset\E$ of an end we define the capacity (\cite{Grig}, pg 152)
$$
\cpy_{\E} \left( \F \right) = \inf \left\{ \int_{\E} \left| \nabla u \right|^2 \, d\mu \right\} 
$$
where the infimum is over locally Lipschitz functions $u$ satisfying $\left. u \right|_{\F} = 1$ and $u$ going to zero at infinity. Then the Maz'ja constant of $\E$ is 
$$
m \left( \E \right) = \inf_{\F\subset\subset\E} \frac{\cpy_{\E} (\F)}{\mu (\F)} \, ,
$$
where $\mu (\F)$ is the volume of $\F$. Maz'ja showed \cite{Maz,Grig2} that there is a constant $c$ such that
$$
c\cdot m\left( \E \right) \le \lambda_0 \left( \E \right) \le m \left( \E \right) 
$$
so that by theorem \ref{Bai1} discreteness is equivalent to the unboundedness of the Maz'ja constant down the ends. 
On the other hand if we consider, for instance in the infinite volume case, the function $v(r;t,s)$, $t<s$,
$$
v (r) = \left\{ \begin{array}{cl}
1 & : t<r<s \\
\frac{\int^{\infty}_r \bar{\m}^{-1} \, dr'}{\int^{\infty}_s \bar{\m}^{-1} \, dr'} & : s<r
\end{array} \right.
$$
in the definition of the capacity we see that the capacity of the subset $\E_t\setminus\E_s$ is bounded above by
$$
\left( \int^{\infty}_s \bar{\m}^{-1} \, dr' \right)^{-1} \, .
$$
The volume of $\E_t\setminus\E_s$ is $\int^{s}_t \bar{\m} \, dr$ so that the terms which appear in (\ref{ldc2}) are in fact an estimate of the Maz'ja constant. Maz'ja's general condition for discreteness then becomes our simpler condition for discreteness subject to the hypotheses \ref{hyp1} and \ref{hyp2} on the structure of our manifold.
\section{Brownian Motion on Manifolds with Ends}
A discussion of brownian motion on manifolds may be found in \cite{Grig}. Here we discuss the relationship between stochastic (in)completeness and discreteness (for a precise definition see \cite{Grig}, a rough definition is that a manifold is stochastically incomplete if brownian motion on the manifold escapes to infinity in a finite time). 
\begin{thm} 
If a manifold $\M$, satisfying hypotheses \ref{hyp1} and \ref{hyp2}, with only one end has essential spectrum then $\M$ is stochastically complete.
\end{thm}
\begin{proof}
Using lemma \ref{pertlem} we see that the perturbed metric (\ref{pertm}) also has a Laplacian with essential spectrum. We choose the constant $c$ in (\ref{pertm}) to be the same as the constant appearing in hypothesis \ref{hyp2} and denote the square root of the determinant of the perturbed metric by $\bar{\m}_c = e^{cr} \bar{\m}$. \\
From lemma \ref{Bai2} and theorem \ref{mthm} we have the existence of a $c_0>0$ such that
\begin{eqnarray*}
c_0 & = & \lim_{t\to\infty} \left( \lambda_0 \left( \E_t \right) \right)^{-1} \le \underline{\lim}_{t\to\infty} \left( \mu_0 \left( \E_t \right) \right)^{-1} \\
& \le & 4 \, \underline{\lim}_{t\to\infty} \sup_{s>t} \int^{\infty}_s \bar{\m}^{-1}_c \, dr \int^{s}_t \bar{\m}_c \, dr \, .
\end{eqnarray*}
That is we can find $T$ such that for all $t_0 >T$
$$
\frac{c_0}{8} < \sup_{s>t_0} \int^{\infty}_s \bar{\m}^{-1}_c \, dr \int^{s}_{t_0} \bar{\m}_c \, dr \, .
$$
This in turn means that we can find an $s_0>t_0=T$ such that 
$$
0 < c_1 < \int^{\infty}_{s_0} \bar{\m}^{-1}_c \, dr \int^{s_0}_{t_0} \bar{\m}_c \, dr \, .
$$
Putting $t_1=s_0$ and repeating we see that this implies that we can find an increasing sequence $\left\{ t_n \right\}$ such that
$$
0 < c_1 < \int^{\infty}_{t_{n+1}} \bar{\m}^{-1}_c \, dr \int^{t_{n+1}}_{t_n} \bar{\m}_c \, dr =
 \int^{\infty}_{t_{n+1}} \bar{\m}^{-1}_c (r) \int^{t_{n+1}}_{t_n} \bar{\m}_c (s) \, ds \, dr \, .
$$
The sum over this sequence will then diverge
\begin{eqnarray}
\infty & = & \lim_{N\to\infty} \sum^N_{n=0} \int^{\infty}_{t_{n+1}} \bar{\m}^{-1}_c (r) \int^{t_{n+1}}_{t_n} \bar{\m}_c (s) \, ds \, dr \nonumber \\
& \le & \int^{\infty}_{0} \bar{\m}^{-1}_c (r) \int^{r}_{0} \bar{\m}_c (s) \, ds \, dr \, . \label{div}
\end{eqnarray}
Now we consider the solution $u(r)$ of
$$
\mbox{} - u'' - \frac{\bar{\m}'_c}{\bar{\m}_c} u' + u = 0
$$
with boundary values $\left. u \right|_0 = 1$, $\left. u' \right|_0 = 0$. We see that $u$ satisfies
$$
u u' \bar{\m}_c = \int^r_0 \left( u^2 + \left( u' \right)^2\right) \bar{\m}_c \, dr'
$$
so that $u,u'>0$. The solution can also be written 
$$
u(t) = 1 + \int^t_0 \bar{\m}^{-1}_c (r) \int^r_0 u(s) \bar{\m}_c (s) \, ds \, dr
$$
which along with (\ref{div}) implies that $\lim_{r\to\infty} u(r) = \infty$. \\
Then from hypothesis \ref{hyp2} it is clear that $u$ is 1-superharmonic with respect to the (unperturbed) Laplacian
\begin{eqnarray*}
\left(\LB + 1 \right) u & = & \mbox{} - u'' - \frac{\m'}{\m} u' + u \\
& \ge & \mbox{} - u'' - \left( \frac{\bar{\m}'}{\bar{\m}} + c \right) u' + u = \mbox{} - u'' - \frac{\bar{\m}'_c}{\bar{\m}_c} u' + u = 0 \, .
\end{eqnarray*}
Finally, the existence of a divergent 1-harmonic function is equivalent to stochastic completeness (corollary 6.6 of \cite{Grig}) completing the theorem.
\end{proof}
This theorem can be generalised to a manifold with many ends, one has to specify that the Laplacian restricted to each end must have essential spectrum. Conversely we may say that a manifold which is stochastically incomplete necessarily has an end such that the Laplacian restricted to that end is discrete (and in fact the discrete end will have infinite volume. The end structure, ie. hypotheses \ref{hyp1} and \ref{hyp2}, is essential).

\end{document}